\documentclass[12pt,a4paper]{article}

\usepackage{amsmath,amsfonts}
\usepackage{mathrsfs}

\newtheorem{theorem}{Theorem}

\let\goth\mathfrak

\def \Limsup{\mathop{\overline{\lim}}\limits}
\def \Liminf{\mathop{\underline{\lim}}\limits}
\def\Pb{\mathbf{P}} 

\def\Ex{\mathbf{E}}

\def\DD{\mathbb{D}}
\def\SS{\mathbb{S}}

\def\RR{\mathbb{R}}

\def\MM{\mathbb{M}}
\def\II{\mathbb{I}}
\def\JJ{\mathbb{J}}

\def\sgn{{\rm sgn}} 
\def\1{\mbox{1\hspace{-.25em}I}}

\begin{document}
\title{On Multi-Step MLE-Process for Ergodic Diffusion}
\author{Yu.A. Kutoyants\\
{\small Laboratoire de Statistique et Processus, Universit\'e du Maine}\\
{\small  Le Mans,  France and}\\
{\small  International Laboratory of Quantitive Finance, Higher School of Economics}\\
{\small   Moscow, Russia and}\\
{\small  National Research University ``MPEI''}\\
{ \small  Moscow, Russia}\\ 
}

\date{}

\maketitle
%{Running title : }
\begin{abstract}
We propose a new method of the construction of the asymptotically efficient
estimator-processes asymptotically equivalent to the MLE and the same time much
more easy to calculate. We suppose that the observed process is ergodic
diffusion and that there is a learning time interval of the length negligeable
with respect to the whole time of observations. The preliminary estimator
obtained after the learning time is then used in the construction of one-step
and two-step MLE processes. We discuss the possibility of the applications of
the proposed estimation procedure  to several other observations models. 

\end{abstract}
%\noindent MSC 2000 Classification: 62M02,  62G10, 62G20.


\section{Introduction}

We consider the problem of parameter  estimation by the continuous time  observations
$X^T=\left(X_t,0\leq t\leq T\right)$ of the
diffusion process 
\begin{align*}
{\rm d}X_t=S\left(\vartheta ,X_t\right)\,{\rm d}t+\sigma
\left(X_t\right)\,{\rm d}W_t,\qquad X_0,\quad 0\leq t\leq T. 
\end{align*}
We suppose that the process $\left(X_t\right)_{t\geq 0}$ has ergodic
properties  with the
density of invariant distribution $f\left(\vartheta ,x\right)$. The functions
$S\left(\vartheta ,x\right)$ and $\sigma \left(x\right)$ are known, smooth and the
parameter $\vartheta \in \Theta\subset{\cal  R}^d $. It is known that the
maximum likelihood estimator (MLE) $\hat\vartheta _T$ constructed by the
observations $X^T$ (under regularity conditions) is
consistent, asymptotically normal and asymptotically efficient (see, e.g.,
\cite{Kut77} or \cite{Kut04}).

We consider here slightly different statement of the problem. Suppose that we
are interested by an {\it estimator-process} $\bar\vartheta
^T=\left(\bar\vartheta _{t,T}, 0\leq t\leq T\right)$, where $\bar\vartheta
_{t,T}$ depends on the observations $X^t=\left(X_s,0\leq s\leq t\right)$ only.
The need of such {\it on-line} estimators naturally arises in many problems,
for example in  adaptive
control. We used such estimators in the construction of the approximation of
the solution of backward stochastic differential equation \cite{KZ14},
\cite{Kut14}, where  the estimator-processes were of the one-step MLE-type. Note that
there is a large literature on stochastic approximation, which provides
satisfactory solutions (see, e.g.;  \cite{KY},  and references therein). For
continuous time systems such problems were studied for example in \cite{H67}
and \cite{NH73}. In the last work there was proposed a recurrent
asymptotically efficient estimation in the case of observations
\begin{align*}
{\rm d}X_t=S\left(\vartheta ,t\right){\rm d}t+\sigma \left(t\right){\rm
  d}W_t,\quad 0\leq t\leq T. 
\end{align*} 
Another recurrent estimator-process for diffusion processes was studied in
\cite{LSZ94}. 
We have to mention that the estimator-process proposed in our work is not
recurrent in this sense. The right-hand side of the equation for it does not depend
on the preceding values of this estimator (see \eqref{01} below).

  For diffusion processes observed in discrete times the
adaptive estimation of the parameters of the trend  and diffusion coefficients
were studied in many works, see, e.g., \cite{Y92}, \cite{Kes97}, \cite{UY12},\cite{UY14}
and the references therein.  Note that in the works \cite{UY12} and
\cite{UY14} the proposed multi-step adaptive procedures of parameter
estimation allows to improve the initial {\it bad} rates of convergence up to
asymptotically efficient ({\it good})  rates.   

The studied in the present work estimator-processes are based on the one-step
MLE structure.  Recall that the one-step MLE was introduced by Le Cam \cite{LC56} in
1956. The definition and properties of it in i.i.d. case can be
found, for example,  in \cite{LR05}. Let us remind it's construction. Suppose
that the observed i.i.d. random variables $X^n=\left(X_1,\ldots,X_n\right)$
have smooth  density function $f\left(\vartheta ,x\right)$ with unknown parameter
$\vartheta \in \Theta \subset {\cal R}^d$ and we have to estimate $\vartheta
$. Suppose as well that we are given an
estimator $\bar\vartheta _n$ which is consistent and asymptotically normal
$\sqrt{n} \left(\bar\vartheta _n-\vartheta _0\right)\Rightarrow {\cal
  N}\left(0,\DD\left(\vartheta _0\right)\right)$ with a {\it bad} $d\times d$ limit
covariance matrix $\DD\left(\vartheta _0\right)>\II\left(\vartheta
_0\right)^{-1}$, i.e., the matrix $\DD\left(\vartheta _0\right)-\II\left(\vartheta
_0\right)^{-1} $ is positive definite. We say ``bad'' because
$\DD\left(\vartheta _0\right)$ is not equal to the the  inverse Fisher information matrix
 $\II\left(\vartheta
_0\right)^{-1} $, which is limit covariance of asymptotically efficient
estimators. The one-step MLE is defined as follows
\begin{align*}
\vartheta _n^\star=\bar\vartheta _n+n^{-1}\II\left(\bar\vartheta _n
\right)^{-1}\sum_{j=1}^{n}\frac{\dot f\left(\bar\vartheta _n
  ,X_j\right)}{f\left(\bar\vartheta _n ,X_j\right)}.
\end{align*}
 Here and in the sequel dot means derivation w.r.t. $\vartheta $. The
 estimator $\vartheta _n^\star$ has 
 already the {\it good} limit covariance matrix
\begin{align*}
%\label{0-1}
\sqrt{n}\left(\vartheta _n^\star-\vartheta _0\right)\Longrightarrow {\cal
  N}\left( 0, \II\left(\vartheta _0\right)^{-1}\right).
\end{align*}
Therefore this procedure allows us to improve any estimator with {\it good}  rate
$\sqrt{n}$ but {\it bad}  limit covariance up to asymptotically efficient.

The one-step MLE $\vartheta _T^\star$ for ergodic diffusion process  can be
defined by a similar way. Suppose that we have a preliminary estimator
$\bar\vartheta _T$ (say, 
minimum distance estimator or estimator of the method of moments),
which is consistent and asymptotically normal : 
$$ 
\sqrt{T}\left( \bar\vartheta
_T-\vartheta \right) \Rightarrow {\cal N}\left(0,\DD\left(\vartheta
\right)\right).
$$
The limit variance $\DD\left(\vartheta
\right)>\II\left(\vartheta \right)^{-1}$, where $\II\left(\vartheta \right) $
is the Fisher information matrix
\begin{align*}
\II\left(\vartheta \right)=\int_{-\infty }^{\infty } \frac{\dot
  S\left(\vartheta ,x\right)\dot S\left(\vartheta ,x\right)^* }{\sigma
  \left(x\right)^2}f\left(\vartheta ,x\right) \;{\rm d}x.
\end{align*}
Here and in the sequel $A^*$ means transposition of $A$. 

 Following the same ``one-step''  idea we can  improve this
estimator up to asymptotically efficient as follows
\begin{align}
\label{01}
\vartheta _T^\star=\bar\vartheta_T+T^{-1/2} \II\left(\bar\vartheta_T
\right)^{-1}\int_{0}^{T} \frac{\dot
  S\left(\bar\vartheta_T,x\right) }{\sigma
  \left(x\right)^2} \;\left[{\rm d}X_t- S\left(\bar\vartheta_T,x\right){\rm
    d}t\right] .
\end{align} 
whehe $\vartheta _T^\star $ is the one-step MLE.
Of course, the special attention have to be paid for the definition of the
stochastic integral because the estimator $\bar\vartheta_T $ depends on the
whole observations $X^T$. This estimator is consistent, asymptotically  
normal 
\begin{align*}
\sqrt{T}\left( \vartheta
_T^\star-\vartheta \right) \Rightarrow {\cal N}\left(0,\II\left(\vartheta
\right) ^{-1}\right)
\end{align*}
 and asymptotically efficient (see, e.g., \cite{Kut04}). Note that the
 preliminary estimator here has a {\it  good} $\sqrt{T}$ rate of convergence.

Recently Kamatani and Uchida \cite{KU14} considered the problem of parameter
estimation by the discrete time observations $X^n=\left(X_{t_i^n},t_i^n=ih_n,
0\leq i\leq n\right)$ of ergodic diffusion process 
\begin{align*}
{\rm d}X_t=S\left(\vartheta_1 ,X_t\right)\,{\rm d}t+\sigma
\left(\vartheta_2 ,X_t\right)\,{\rm d}W_t,\qquad X_0,\quad 0\leq t\leq T. 
\end{align*}
Here $\vartheta =\left(\vartheta _1,\vartheta _2\right)$ is unknown
parameter. They proposed a modification of Newton-Raphson (N-R) procedure,
with the initial estimators of {\it bad} (non-optimal) rate of convergence and
showed that the multi-step N-R procedure allows to obtain the asymptotically
efficient estimators with the {\it good} rates. The asymptotic is
$h_n\rightarrow 0$ and $nh_n\rightarrow \infty $.

In our work we consider a similar  construction but  based on the modification of
one-step MLE procedure. We propose    estimator-processes  $\vartheta
_{t,T}^\star,T^\delta \leq t\leq T $ and  $\vartheta
_{t,T}^{\star\star},T^\delta \leq t\leq T $, where $\delta <1$ and $\vartheta _{t,T}^\star $
and $\vartheta
_{t,T}^{\star\star} $ have one-step MLE-type structure. As
preliminary estimator $\bar\vartheta _{T^\delta } $  we take an estimator
constructed by the first observations $X^{T^\delta }$  on the time interval $\left[0,T^\delta
  \right]$. Therefore the preliminary estimator has a {\it bad} rate of
convergence due to the length of the learning interval. 

  Then we propose one-step (for $\delta \in
\left(\frac{1}{2},1\right)$) and two-step (for $\delta \in
(\frac{1}{4},\frac{1}{2}]$) MLE-processes.  
 For example, the one-step MLE-process is
\begin{align*}
\vartheta _{t,T}^\star=\bar\vartheta _{T^\delta } +T^{-1/2}\II\left(\bar\vartheta _{T^\delta }
\right) ^{-1/2}\int_{T^\delta }^{t} \frac{\dot
  S\left(\bar\vartheta_{T^\delta },x\right) }{\sigma
  \left(x\right)^2} \;\left[{\rm d}X_t- S\left(\bar\vartheta_{T^\delta },x\right){\rm
    d}t\right],
\end{align*}
where $ T^\delta \leq t\leq T $. It is shown that this estimator-process is
consistent, asymptotically normal and asymptotically efficient. Note that the
calculation of this estimator-process is much more simple than the calculation of
the MLE $\hat\vartheta _{t,T}$ for all $t\in [T^\delta ,T]$.

\section{Auxiliary results}

We are given a probability space $\left\{\Omega , {\cal F}, \Pb\right\}$ with
filtration $\left({\cal F}_t\right)_{t\geq 0}$ satisfying the {\it usual
  conditions} and the Wiener process $W=\left(W_t,{\cal F}_t,t\geq
0\right)$. Suppose that for all $\vartheta \in\Theta\subset \RR^d $, ($\Theta
$ is an open bounded set) the stochastic process $X=\left(X_t,{\cal F}_t,t\geq
0\right)$ satisfies the stochastic differential equation
\begin{equation}
\label{2-1}
{\rm d}X_t=S\left(\vartheta ,X_t\right)\,{\rm d}t+\sigma
\left(X_t\right)\,{\rm d}W_t,\quad X_0,\quad  t\geq 0,
\end{equation}
where $\left(X_0,{\cal F}_0 \right)$ is the initial value. The functions $S\left(\vartheta
,x\right)$ and $\sigma \left(x\right)^2>0$ are such 
that this equation has a unique strong solution on any fixed interval
$\left[0,T\right]$ and that the measures $\left\{\Pb_\vartheta
^{\left(T\right)},\vartheta \in \Theta \right\}$ induced in the measurable
space $\left({\cal C}\left[0,T\right],{\goth B}_T\right)$ of its realizations
are equivalent (see conditions for example here \cite{LS05}).  Moreover we
suppose that  the process
$X_t,t\geq 0$ has ergodic properties  with the density of invariant
distribution 
\begin{align*}
f\left(\vartheta ,x\right)=\frac{1}{G\left(\vartheta \right)\sigma
  \left(x\right)^2}\exp\left\{2 \int_{0}^{x} \frac{S\left(\vartheta
  ,y\right)}{\sigma \left(y\right)^2 }\;{\rm d}y\right\},
\end{align*}
where $G\left(\vartheta \right)$ is the normalizing constant \cite{Kh12}. The
random variable with such density we denote as $\xi $ and suppose that $X_0$
has the same density function. This condition makes the stochastic process
stationary. 

The sufficient condition for the existence of ergodic properties 
 we take as in  \cite{Kut04}. Define the class of functions
\begin{align*}
{\cal P}=\left\{h\left(\cdot \right):\quad \left|h{\left(x\right)}\right|\leq
C\left(1+\left|x\right|^q\right)\right\}, 
\end{align*}
where the constants $C>0,q>0$  do
not depend on $\vartheta $ in the case of the function $S\left(\vartheta
,x\right)$ and its derivatives and can be different for different functions. 

Condition ${\cal A}_0\left(\Theta \right)$.  {\it The functions $S\left(\vartheta
,\cdot \right) , \sigma \left(\cdot \right)^{\pm 1} \in {\cal P}$  and}
\begin{align*}
\Limsup_{\left|x\right|\rightarrow \infty }\sup_{\vartheta \in \Theta
}\frac{\sgn\left(x\right)\,S\left(\vartheta ,x\right)}{\sigma
  \left(x\right)^2}<0.
\end{align*} 
The {\it smoothness}  condition: {\it the function $S\left(\vartheta
  ,x\right)$ has two
continuous partial derivatives w.r.t. $\vartheta $  and these derivatives belong to
${\cal P}$. }

These derivatives we denote as : $\dot S\left(\vartheta
,x\right)$  (vector)   and  $\ddot \SS\left(\vartheta ,x\right)$ ($d\times d$ matrix).

The {\it identifiability} condition: {\it for any} $\nu >0$
\begin{align*}
\inf_{\vartheta _0\in \Theta }\inf_{\left|\vartheta -\vartheta _0\right|>\nu
}\Ex_{\vartheta _0}\left(\frac{S\left(\vartheta ,\xi_0
  \right)-S\left(\vartheta_0 ,\xi_0 \right)}{\sigma \left(\xi_0
  \right)}\right)^2>0.
\end{align*}
Here the r.v. $\xi_0$ has the
density function $f\left(\vartheta _0,x\right)$. 

The Fisher information matrix ${\II} \left(\vartheta\right) $
 is uniformly non
degenerate  (below $\lambda \in R^d$) 
\begin{align*}
\inf_{\vartheta \in \Theta }\inf_{\left|\lambda \right|=1 }\lambda ^*{\II}
\left(\vartheta \right)\lambda  >0.
\end{align*}

The set of all these conditions we call {\it Regularity conditions}.

We have to estimate $\vartheta $ by the observations $X^t=\left(X_s,0\leq s\leq
t\right)$ for  all $ t\in(0,T]$ and to describe the properties of the
  estimator-process $ \left(\bar\vartheta\left(t\right),{\cal F}_t,\right.$
$\left.0< t\leq
  T\right)$, where 
$\bar\vartheta\left(t\right)=\bar\vartheta
\left(t,X^t\right)$.  We would like to
obtain an estimator-process which has asymptotically optimal in some sense
properties.

 It will be convenient to change the variables $t=\tau T, \tau
\in(0,1]$ and to study the random processes $\bar\vartheta
_{\tau,T},0< \tau \leq 1$, where $\bar\vartheta _{\tau,T}= \bar\vartheta \left(\tau
T\right)$. For simplicity of exposition we will write $ \bar\vartheta
_{\tau,T}  $ as $\bar\vartheta _{\tau}$.

Introduce the likelihood ratio-process $V\left(\vartheta ,X^{\tau
  T}\right),\vartheta \in \Theta,\;0<\tau \leq 1  $, where
\begin{align*}
V\left(\vartheta ,X^{\tau T}\right)=\exp\left\{\int_{0}^{{\tau T}}\frac{S\left(\vartheta
  ,X_s\right)}{\sigma \left(X_s\right)^2} \,{\rm d}X_s-\int_{0}^{{\tau T}}\frac{S\left(\vartheta
  ,X_s\right)^2}{2\sigma \left(X_s\right)^2} \,{\rm d}s\right\} . 
\end{align*}

 Note that for any $\tau \in (0,1]$ the family of measures
   $\left\{\Pb_\vartheta ^{\left(\tau T\right)},\vartheta \in \Theta \right\}$
   is {\it locally asymptotically normal} (LAN) in $\Theta $, i.e.; the
   likelihood ratio-process 
\begin{align*}
Z_{\tau T}\left(u\right)=\frac{V\left(\vartheta _0+\frac{u}{\sqrt{T}}, X^{\tau
      T}\right)}{V\left(\vartheta _0, X^{\tau T}\right) },\qquad 0<\tau \leq 1,
\end{align*}
with $u$ such that $\vartheta _0+\frac{u}{\sqrt{T}}\in\Theta  $, admits the representation 
\begin{align*}
Z_\tau \left(u\right)=\exp\left\{ u^*\tilde\Delta _\tau \left(\vartheta _0,X^{\tau
  T}\right)-\frac{1}{2} u^*\II_\tau \left(\vartheta _0\right)u+r_T\right\}.
\end{align*}
Here  $\II_\tau \left(\vartheta _0\right)=\tau\II \left(\vartheta _0\right) $, $r_T\rightarrow 0$ and the score-function (vector-process)
\begin{align*}
\tilde\Delta _\tau \left(\vartheta _0,X^{\tau
  T}\right)=\frac{1}{\sqrt{T}}\int_{0}^{\tau T}\frac{\dot S\left(\vartheta
  _0,X_s\right)}{\sigma \left(X_s\right)^2} \left[{\rm d}X_s-S\left(\vartheta
  _0,X_s\right){\rm d}s\right] \Longrightarrow {\cal N}\left(0,\II_\tau
\left(\vartheta _0\right)\right).
\end{align*}

  Therefore we
   have the following Hajek-Le Cam-type low bound for polynomial loss function
   ($p>0$): for all estimator-processes $ \bar\vartheta _{\tau},0<\tau \leq 1$
   and all $\vartheta _0\in \Theta $ and $\tau \in(0,1]$ 
\begin{align}
\label{2-2}
\Liminf_{\nu \rightarrow 0}\Liminf_{T\rightarrow \infty }\sup_{\left|\vartheta
  -\vartheta _0\right|<\nu } \Ex_\vartheta \left|\sqrt{ T}{\II}_\tau
\left(\vartheta_0 \right)^{1/2}\left(\bar\vartheta _\tau -\vartheta
\right)\right|^p\geq \Ex\left|\zeta \right|^p,
\end{align}
where  the vector  $\zeta \sim {\cal N}\left(0,\JJ\right)$, $\JJ$ is unit $d\times
d$ matrix (see, e.g., \cite{IH81} or  \cite{Kut04}).

Therefore the estimator process $\bar\vartheta _\tau,0<\tau \leq 1$ we call
asymptotically efficient if for all $\vartheta _0\in \Theta $ and all $\tau
\in (0,1]$ we have the equality    
\begin{align}
\label{2-3}
\lim_{\delta \rightarrow 0}\lim_{T\rightarrow \infty
}\sup_{\left|\vartheta -\vartheta _0\right|<\delta }\Ex_\vartheta
\left|\sqrt{T}{\II}_\tau \left(\vartheta_0\right)^{1/2}\left(\bar\vartheta _\tau -\vartheta
\right)\right|^p=\Ex\left|\zeta \right|^p.
\end{align}
 One solution of this problem is to introduce the MLE process $\hat\vartheta
 _\tau, 0<\tau \leq 1 $ defined by the equation
\begin{align}
\label{2-4}
V\left(\hat\vartheta _\tau ,X^{\tau T}\right)=\sup_{\vartheta \in \Theta
}V\left(\vartheta ,X^{\tau T}\right), \qquad \tau \in(0,1].
\end{align}

It is known that the estimators $\hat\vartheta _\tau ,\tau \in(0,1]$ under regularity
conditions are consistent and asymptotically normal ($\vartheta _0$ is the
true value)
\begin{align*}
\sqrt{T}\left(\hat\vartheta _\tau -\vartheta _0\right)\Longrightarrow {\cal
  N}\left(0,{\II}_\tau \left(\vartheta _0\right)^{-1}\right).
\end{align*}
Moreover we have the uniform on  $\vartheta $ convergence of moments
$$ \lim_{T\rightarrow \infty }T^{p/2}\Ex_\vartheta \left|{\II}_\tau
\left(\vartheta \right)^{1/2}\left(\hat\vartheta _\tau -\vartheta
\right)\right|^p=\Ex\left|\zeta \right| ^p.
$$ 
Therefore the MLE-process $\hat\vartheta _\tau ,0<\tau \leq 1 $ is
asymptotically efficient (see \cite{Kut04}).

 Unfortunately except the linear case $S\left(\vartheta ,x\right)=\vartheta^*
 h\left(x\right)$ the calculation for all $\tau \in(0,1]$ of the MLE-process
   as solution of the equation \eqref{2-3} is computationally very difficult
   problem and we have to seek another estimator-process which is
   computantionally much more simple. The goal of this work is to describe such
   class of estimator-processes. The construction of the proposed
   estimator-processes is based on the development of the well-known one-step
   MLE device.

\section{Main result}

We consider the problem of estimation $\vartheta $ by observations $X^t$ for
$t\in \left[T^\delta ,T\right]$. The corresponding estimators we study after
the change of variables $t=\tau T$. Therefore we are interested by the
construction of the estimator-process $\vartheta _{\tau ,T}^\star, \tau \in
\left[\tau _\delta ,1\right]$, where $\tau _\delta =T^{-1+\delta }$. We show
that if $\delta \in (\frac{1}{2},1)$, then the one-step MLE-process is
  asymptotically normal and asymptotically efficient. If $\delta \in
  (\frac{1}{4},\frac{1}{2}]$, then we propose two-step MLE-process with the same
    asymptotic properties. 

\subsection{One-step MLE ($\delta \in (\frac{1}{2},1)$)}

Introduce the learning interval $0\leq t\leq T^\delta $, where $\delta \in
(\frac{1}{2},1)$  and denote by $\bar \vartheta _{\tau _\delta } $ an estimator of
parameter $\vartheta $ which  is uniformly  on compacts
${\bf A}\subset\Theta $ asymptotically normal 
 \begin{equation}
\label{2-5}
T^{\frac{\delta }{2}}\left(\bar \vartheta _{\tau_\delta } -\vartheta _0
\right)\Longrightarrow {\cal N}\left(0,\DD\left(\vartheta _0\right)\right),
\end{equation}
where $\tau _\delta =T^{-1+\delta }\rightarrow 0$ and the matrix
$\DD\left(\vartheta _0\right) $ of limit covariance is 
bounded. Moreover we suppose that we have the convergence of all polynomial
moments too: for all $p>0$ 
\begin{equation}
\label{2-6}
\sup_{\vartheta_0 \in {\bf A}} T^{\frac{p\delta }{2}}\Ex_{\vartheta_0} \left|\bar
\vartheta _{\tau_\delta } -\vartheta_0 \right|^p<C ,
\end{equation}
where the constant  $C>0$ does not depend on $T$. The regularity conditions
providing these properties of the MLE, minimum distance estimators (MDE), bayesian
estimators (BE) and the estimators of the method of moments (EMM) can be found, for
example,  in
\cite{Kut04}, Chapter 2. Therefore as preliminary estimator we can take one of
them.

The one-step MLE-process we construct as follows:
\begin{equation}
\label{2-7}
\vartheta _\tau ^\star=\bar\vartheta _{\tau_\delta }+\frac{\II\left(
  \bar\vartheta _{\tau_\delta }\right)^{-1}}{ \sqrt{\tau  T}  }\Delta _\tau 
\left(\bar\vartheta _{\tau_\delta },X^{\tau T}_{T^\delta
} \right),\quad \tau \in
\left[\tau _\delta ,1\right] 
\end{equation}
where
\begin{align}
\label{2-8}
\Delta _\tau
\left(\vartheta ,X^{\tau T}_{T^\delta
} \right)&=\frac{1}{\sqrt{{\tau  T}}}\int_{T^\delta
}^{\tau T}\frac{\dot S\left(\vartheta ,X_t\right)}{\sigma
  \left(X_t\right)^2}\left[{\rm d}X_t-S\left(\vartheta ,X_t\right){\rm
    d}t\right].
\end{align}

Introduce the random process 
$$
\eta _{\tau ,T}\left(\vartheta _0\right) = \tau \sqrt{ T} \II\left(\vartheta
_0\right)^{1/2}\left( \vartheta_\tau 
^\star-\vartheta _0\right),\qquad \tau _*\leq \tau \leq 1, 
$$ 
where $\tau _*\in (0,1)$ and measurable space $\left( {\cal C}\left[\tau
  _*,1\right],{\goth B}\right)$ of continuous on $\left[\tau _*,1\right]$
functions. Here ${\goth B} $ is the corresponding borelian $\sigma
$-algebra. Denote by $W\left(\tau \right),0\leq \tau \leq 1$ a
$d$-dimensional standard Wiener process.

\begin{theorem}
\label{T1}
Suppose that  the regularity conditions    hold. Then the one-step MLE-process  $\vartheta
_\tau ^\star, \tau _\delta<\tau \leq 1$ has the following properties:
\begin{enumerate}
\item It  is  uniformly consistent: for any $\nu
>0$
\begin{equation}
\label{2-9}
\lim_{T\rightarrow \infty }\sup_{\vartheta_0 \in {\bf A}}\Pb_{\vartheta _0}\left(\sup_{\tau _\delta
  \leq \tau \leq 1}\left| \vartheta
_\tau ^\star-\vartheta _0\right|>\nu \right)=0.
\end{equation}
\item For any $\tau _*\in (0,1)$ the random process $\eta _{\tau
  ,T}\left(\vartheta _0\right),\tau _*\leq \tau \leq 1 $ converges in
  distribution in $\left( {\cal C}\left[\tau _*,1\right],{\goth B}\right)$ to
  the vector-process $W\left(\tau \right),\tau _*\leq \tau \leq 1$.
\item It is asymptotically efficient in the sense \eqref{2-2} for any $p>0$. 
\end{enumerate}

\end{theorem}

Note that for a fixed $\tau $ we have the asymptotic normality 
\begin{equation}
\label{2-10}
\sqrt{ T}\left( \vartheta
_\tau ^\star-\vartheta _0\right)\Longrightarrow  {\cal N}\left(0,
\II_\tau \left(\vartheta _0\right)^{-1}\right) 
\end{equation}

{\bf Proof.} We denote by $C$ the generic constant. The uniform consistency
\eqref{2-9} is verified as follows
\begin{align*}
&\Pb_{\vartheta _0}\left( \sup_{\tau _\delta \leq \tau \leq 1}\left| \vartheta
  _\tau ^\star-\vartheta _0\right|>\nu \right)\leq \Pb_{\vartheta _0}\left(
  \left| \bar\vartheta _{\tau_\delta } -\vartheta _0\right|>\frac{\nu
  }{2}\right)\\ &\qquad +\Pb_{\vartheta _0}\left( \sup_{\tau _\delta \leq \tau
    \leq 1} \left(\tau  T \right)^{-1/2}\left|\II\left( \bar\vartheta
  _{\tau_\delta }\right)^{-1}\Delta _\tau \left(\bar\vartheta _{\tau_\delta
  },X^{\tau T} \right)\right| >\frac{\nu }{2} \right)\\ 
&\quad \leq \nu
  ^{-p}{2^p\Ex_{\vartheta _0} \left| \bar\vartheta _{\tau_\delta } -\vartheta
    _0\right|^p }{}+\Pb_{\vartheta _0}\left( \sup_{\tau _\delta \leq \tau \leq
    1} \left|\Delta _\tau \left(\bar\vartheta _{\tau_\delta },X^{\tau T}
  \right)\right| >\frac{\nu T^{\delta /2}}{2C} \right)\\ 
&\quad \leq C\nu
  ^{-p}T^{-p\delta/2 }+\Pb_{\vartheta _0}\left(\sup_{\tau _\delta \leq \tau
    \leq 1} \frac{1}{{\sqrt{ T}}} \left|\int_{T^\delta }^{\tau T}
  \frac{\dot S\left(\bar\vartheta _{\tau_\delta },X_t\right)}{\sigma \left(X_t\right)}{\rm
    d}W_t\right|>\frac{\nu T^{\delta /2}}{4C} \right)\\ 
&\qquad
  +\Pb_{\vartheta _0}\left( { }{ }\sup_{\tau _\delta \leq \tau \leq 1} 
 \left| \int_{T^\delta }^{\tau T} \frac{\dot S\left(\bar\vartheta _{\tau_\delta } ,X_t\right)
    \left[S\left(\vartheta _{0} ,X_t\right)-S(\bar\vartheta_{\tau_\delta } ,X_t)\right] }{T\;\sigma
    \left(X_t\right)^2}{\rm d}t\right| >\frac{\nu }{4C} \right).
\end{align*}
We have
\begin{align*}
S\left(\vartheta _{0} ,x\right)-S(\bar\vartheta_{\tau_\delta }
,x)=\int_{0}^{1}\left(\vartheta _{0}- \bar\vartheta_{\tau_\delta
}\right)^*\dot S\left(\vartheta _v,x\right)\, {\rm d}v,
\end{align*}
where $\vartheta _v=\vartheta _0+v\left( \vartheta _{0}- \bar\vartheta_{\tau_\delta
}\right) $
and therefore 
\begin{align}
&\sup_{\tau _\delta \leq \tau \leq 1} 
 \left| \int_{T^\delta }^{\tau T} \frac{\dot S\left(\bar\vartheta _{\tau_\delta } ,X_t\right)
    \left[S\left(\vartheta _{0} ,X_t\right)-S(\bar\vartheta_{\tau_\delta } ,X_t)\right] }{T\;\sigma
    \left(X_t\right)^2}{\rm d}t\right|\nonumber\\
&\qquad \leq \int_{T^\delta }^{ T} \frac{\left|\dot S\left(\bar\vartheta
   _{\tau_\delta } ,X_t\right)\right|\; 
    \left|S\left(\vartheta _{0} ,X_t\right)-S(\bar\vartheta_{\tau_\delta } ,X_t)\right| }{T\;\sigma
    \left(X_t\right)^2}{\rm d}t\nonumber\\
&\qquad \leq \int_{0}^{1}  \int_{T^\delta }^{ T} \frac{\left|\dot
   S\left(\bar\vartheta _{\tau_\delta } ,X_t\right)\right|\; 
    \left|\dot S\left(\vartheta _v ,X_t\right)\right| }{T\;\sigma
    \left(X_t\right)^2}{\rm d}t\,{\rm d}v \left|\vartheta _{0}-
 \bar\vartheta_{\tau_\delta} \right|\nonumber\\ 
&\qquad \leq \frac{C}{T}\int_{T^\delta }^{
   T}\left(1+\left|X_t\right|^q\right){\rm d}t\;\left|\vartheta _{0}-
 \bar\vartheta_{\tau_\delta} \right|\longrightarrow 0,
\label{2-11}
\end{align}
where we used the condition $\dot S\left(\vartheta ,x\right),\sigma
\left(x\right)^{-1}\in {\cal P}$ and the consistency of $\bar\vartheta_{\tau_\delta}  $.
 Recall that by condition ${\cal A}_0\left(\Theta \right)$ the invariant
density has all polynomial moments and therefore we obtain the convergence to
zero of
all moments for these normalized  integrals.

Further, we have for any $\lambda \in R^d$ and $m>0$
\begin{align*}
&\Pb_{\vartheta _0}\left(\sup_{\tau _\delta \leq \tau \leq 1} \frac{1}{{\sqrt{
        T}}} \left|\int_{T^\delta }^{\tau T} \frac{\lambda^*\dot
    S\left(\bar\vartheta_{\tau_\delta} ,X_t\right)}{\sigma
    \left(X_t\right)}{\rm d}W_t\right|>\frac{\nu T^{\delta /2}}{4C}
  \right)\\
 &\qquad \leq \Pb_{\vartheta _0}\left(\sup_{\tau _\delta \leq \tau
    \leq 1} \left|\int_{T^\delta }^{\tau T} \frac{\lambda^*\dot
    S\left(\bar\vartheta_{\tau_\delta} ,X_t\right)}{\sigma
    \left(X_t\right)}{\rm d}W_t\right|>\frac{\nu T^{\frac{\delta+1}{2} }}{4C}
  \right)\\ 
&\qquad \leq C\nu ^{-2m}T^{-m\left(\delta+1\right) }\Ex_{\vartheta
    _0}\left|\int_{T^\delta }^{ T} \frac{\lambda^*\dot
    S\left(\bar\vartheta_{\tau_\delta} ,X_t\right)}{\sigma
    \left(X_t\right)}{\rm d}W_t\right|^{2m}\\
 &\qquad \leq C\nu
  ^{-2m}T^{-m\left(\delta+1\right) }\Ex_{\vartheta _0}\left(\int_{T^\delta }^{
    T} \left|\frac{\lambda^*\dot S\left(\bar\vartheta_{\tau_\delta}
    ,X_t\right)}{\sigma \left(X_t\right)}\right|^2{\rm
    d}t\right)^{m}\leq \frac{C}{\nu ^{2m}T^{m\delta}}\rightarrow 0,
\end{align*}
where we used the  Burkholder-Davis-Gundy (BDG) inequality (see, e.g.,
\cite{KS}, Theorem 3.28). 
 Therefore the consistency \eqref{2-8} is proved.

To prove the weak convergence of $\eta _{\tau ,T}\left(\vartheta
_0\right),\tau _*\leq \tau \leq 1 $ we write the representation 
\begin{align*}
\eta _{\tau ,T}\left(\vartheta _0\right)=\tilde\eta _{\tau ,T}\left(\vartheta
_0\right)+P\left(\tau ,X^{\left(\tau T\right)}\right)+R\left(\tau
,X^{\left(\tau T\right)}\right), 
\end{align*}
where
\begin{align*}
 \tilde \eta _{\tau ,T}\left(\vartheta _0\right)=\II\left(\vartheta
 _0\right)^{-1/2}\frac{1}{\sqrt{T}} \int_{T^\delta }^{\tau T}\frac{\dot 
    S\left(\bar\vartheta_{T^\delta } ,X_t\right)}{\sigma    \left(X_t\right)}\,{\rm d}W_t.
\end{align*}
Then we verify that 
\begin{enumerate}
\item We have the uniform convergence
\begin{align}
\label{2-12}
\sup_{\tau _*\leq \tau \leq 1}\left|P\left(\tau ,X^{\left(\tau
  T\right)}\right) \right|\rightarrow 0,\qquad \sup_{\tau _*\leq \tau \leq
  1}\left|R\left(\tau ,X^{\left(\tau T\right)}\right) \right|\rightarrow 0. 
\end{align}
\item We   have 
 the convergence of finite-dimensional distributions: for any $k=1,2,\ldots $ and $\tau
  _*\leq \tau _1<\ldots <\tau _k\leq 1$
\begin{equation}
\label{2-13}
\Bigl(\tilde\eta _{\tau_1 ,T}\left(\vartheta _0\right),\ldots,\tilde\eta _{\tau_k
  ,T}\left(\vartheta _0\right) \Bigr)\Longrightarrow \Bigl(W\left(\tau
_1\right),\ldots,W\left(\tau _k\right)\Bigr) .
\end{equation}
\item   There exists  a constant $C>0$ such that
\begin{equation}
\label{2-14}
\sup_{\vartheta _0\in {\bf A}} \Ex_{\vartheta _0}\left|\tilde\eta _{\tau_1
  ,T}\left(\vartheta _0\right)-\tilde\eta _{\tau_2 ,T}\left(\vartheta _0\right)
\right|^4\leq C\left|\tau _1-\tau _2\right|^2 .
\end{equation}
\end{enumerate}
 We consider  $T>T_*$ where  $T_*=\tau _*^{\frac{1}{\delta -1}}$, i.e.,  $\tau
_\delta<\tau _*$.  We can write
\begin{align*}
&\tau \sqrt{T} \II\left(\vartheta _0\right)^{1/2} \left( \vartheta _{\tau
  }^\star-\vartheta_0\right)= \tau\sqrt{T} \II\left(\vartheta
  _0\right)^{1/2}\left( \bar\vartheta_{\tau_\delta }-\vartheta
  _0\right)\\ &\qquad + \II\left(\vartheta
  _0\right)^{1/2}\II\left(\bar\vartheta_{\tau_\delta
  }\right)^{-1}\frac{1}{\sqrt{T}} \int_{T^\delta }^{\tau T}\frac{\dot
    S\left(\bar\vartheta_{\tau_\delta } ,X_t\right)}{\sigma
    \left(X_t\right)}\,{\rm d}W_t\\ &\qquad + \II\left(\vartheta
  _0\right)^{1/2}\II\left(\bar\vartheta_{\tau_\delta
  }\right)^{-1}\frac{1}{\sqrt{T}} \int_{T^\delta }^{\tau T}\frac{\dot
    S\left(\bar\vartheta_{\tau_\delta } ,X_t\right)\left[ S\left(\vartheta_{0
      } ,X_t\right)-S\left(\bar\vartheta_{\tau_\delta } ,X_t\right)
      \right]}{\sigma \left(X_t\right)}\,{\rm d}t\\ &\quad =
  \II\left(\vartheta _0\right)^{-1/2}\frac{1}{\sqrt{T}} \int_{T^\delta }^{\tau
    T}\frac{\dot S\left(\vartheta_0 ,X_t\right)}{\sigma
    \left(X_t\right)}\,{\rm d}W_t+\tilde u_{\tau_\delta }T^{-1/2}P\left(\tau ,X^{\tau
    T}\right)+R\left(\tau ,X^{\tau T}\right),
\end{align*}
where we denoted $\tilde u_{\tau_\delta }= \tau \sqrt{T} \II\left(\vartheta
  _0\right)^{1/2}\left( \bar\vartheta_{\tau_\delta }-\vartheta
  _0\right)$ and 
\begin{align*}
P\left(\tau ,X^{{\tau T}}\right)=\sqrt{T}\left[\II\left(\vartheta_0\right)- \frac{1}{\tau
      {T}} \int_{T^\delta }^{\tau T}\frac{\dot 
    S\left(\vartheta_0 ,X_t\right)\dot S\left(\vartheta_0 ,X_t\right)^* }{\sigma
    \left(X_t\right)}\,{\rm d}t\right].
\end{align*}
Here  $R\left(\tau ,X^{\left(\tau T\right)}\right)$ is a difference
of the corresponding  expressions with $\bar\vartheta _{\tau_\delta }$ and $\vartheta _0$. 

Then we remark that by the central limit theorem for ordinary integrals (see
e.g.,\cite{Man68} or  \cite{Kut04}, Proposition 1.23) we have
\begin{align*}
\lambda P\left(\tau ,X^{\tau T}\right)\lambda ^*\Longrightarrow {\cal N}\left(0, D^2\right)
\end{align*}
with some finite  limit variance $D^2$.

Hence
\begin{align*}
&\tilde u_{\tau_\delta }T^{-1/2} P\left(\tau ,X^{\tau T}\right)=\tau
T^{\frac{\delta }{2}} \II\left(\vartheta _0\right)^{1/2}\left( \bar\vartheta_{\tau_\delta
}-\vartheta _0\right)P\left(\tau ,X^{\tau T}\right)T^{-\frac{\delta}{2} }\rightarrow 0
\end{align*}
because $T^{\frac{\delta }{2}}\left( \bar\vartheta_{\tau_\delta
}-\vartheta _0\right)$ and $P\left(\tau ,X^{\tau T}\right)$ are bounded in
probability. To verify that  this  convergence is uniform in $\tau
\in\left[\tau 
  _*,1\right]$ we give some details. Let us denote
\begin{align*}
h\left(\vartheta _0,x\right)&=\lambda \left(\II\left(\vartheta _0\right)-\frac{\dot 
    S\left(\vartheta_0 ,x\right)\dot S\left(\vartheta_0 ,x\right)^* }{\sigma
    \left(x\right)} \right)\lambda ^*\\
H\left(\vartheta _0,x\right)&=\int_{-\infty }^{x}h\left(\vartheta
_0,y\right)f\left(\vartheta _0,y\right){\rm d}y.
\end{align*}
 Here $\lambda \in R^d, \left|\lambda \right|=1$. Then $\Ex_{\vartheta _0}
 h\left(\vartheta _0,\xi_0 \right)=0$ and by It\^o 
 formula  we can write
\begin{align*}
\lambda P\left(\tau ,X^{\tau T}\right)\lambda ^*&=\frac{1}{\tau
  \sqrt{T}}\int_{0}^{\tau T} h\left(\vartheta _0,X_t\right)\,{\rm
  d}t=\frac{1}{\tau
  \sqrt{T}}\int_{X_0}^{X_{\tau T}}h\left(\vartheta _0,x\right){\rm d}x, \\
& \quad -\frac{1}{\tau
  \sqrt{T}}\int_{0}^{\tau T}\frac{2H\left(\vartheta _0,X_t\right)}{\sigma
    \left(X_t\right) f\left(\vartheta 
  _0,X_t\right)} \;{\rm d}W_t .
\end{align*}
Further by BDG inequality: for any $\nu >0$ and $m>0$
\begin{align*}
&\Pb_{\vartheta _0}\left(\sup_{\tau _*\leq \tau \leq 1} \frac{1}{\tau
  \sqrt{T}}\left|\int_{0}^{\tau T}\frac{2H\left(\vartheta _0,X_t\right)}{\sigma
    \left(X_t\right) f\left(\vartheta 
  _0,X_t\right)} \;{\rm d}W_t\right| >\nu T^{\frac{\delta }{2}}  \right)\\
&\qquad \leq K_m\left(\nu\tau _*\right) ^{-2m}T^{-m\delta }\Ex_{\vartheta
    _0}\frac{1}{T^m}\left(\int_{0}^{\tau T}\frac{4H\left(\vartheta _0,X_t\right)^2}{\sigma
    \left(X_t\right)^2 f\left(\vartheta 
  _0,X_t\right)^2} \;{\rm d}t \right)^m\\
&\qquad \leq C\nu^{-2m}T^{-m\delta }\longrightarrow 0.
\end{align*}
The existence of the related moments is verified following the same steps as
in \cite{Kut04}, p.31-33. 
To show that $R\left(\tau ,X^{\left(\tau T\right)}\right)\rightarrow 0 $ we
verify  several estimates like
\begin{align}
&\Ex_{\vartheta _0}\left|\frac{1}{\sqrt{T}} \int_{T^\delta }^{\tau
    T}\frac{\dot S\left(\vartheta_0 ,X_t\right)}{\sigma
    \left(X_t\right)}\,{\rm d}W_t-\frac{1}{\sqrt{T}} \int_{T^\delta }^{\tau
    T}\frac{\dot S\left(\bar\vartheta_{\tau_\delta } ,X_t\right)}{\sigma
    \left(X_t\right)}\,{\rm d}W_t\right|^2 \nonumber\\ &\qquad =\Ex_{\vartheta _0}
  \left[ \Ex_{\vartheta _0} \left. \left|\frac{1}{\sqrt{T}} \int_{T^\delta
    }^{\tau T}\frac{\dot S\left(\vartheta_0 ,X_t\right)
      -\dot S\left(\bar\vartheta_{\tau_\delta } ,X_t\right) }{\sigma
      \left(X_t\right)}\,{\rm d}W_t\right|^2\right|{\cal F}_{T^\delta
    }\right]\nonumber\\ 
&\qquad \leq \frac{1}{{T}}\int_{T^\delta
    }^{\tau T}\Ex_{\vartheta _0} \left|\frac{\dot S\left(\vartheta_0 ,X_t\right)
      -\dot S\left(\bar\vartheta_{\tau_\delta } ,X_t\right) }{\sigma
      \left(X_t\right)} \right|^2{\rm d}t\nonumber\\ 
&\qquad \leq \frac{C}{{T}}\int_{T^\delta
    }^{\tau T}\Ex_{\vartheta _0}
  \left(1+\left|X_t\right|^q\right)\left|\bar\vartheta_{\tau_\delta }-\vartheta_0
  \right|^2\;{\rm d}t=O\left(T^{-{\delta }{}} \right)\rightarrow 0 
\label{2-15}
\end{align}
and (below $ \vartheta _v=\vartheta _0+v\left(\bar\vartheta_{\tau_\delta } -\vartheta _0\right)$)
\begin{align}
&\frac{1}{ {\sqrt{T}}} \left|\int_{0 }^{T^\delta}\frac{\dot
  S\left(\bar\vartheta_{\tau_\delta } ,X_t\right)\left[S\left(\vartheta_{0} ,X_t\right)-
  S\left(\tilde\vartheta_{\tau_\delta } ,X_t\right)\right]}{\sigma
  \left(X_t\right)^2}\,{\rm d}t\right|\nonumber\\
&\qquad =\left|\frac{\tilde u_{\tau_\delta }}{ {T}} \int_{0 }^{1}\int_{0 }^{T^\delta}\frac{\dot
  S\left(\bar\vartheta_{\tau_\delta } ,X_t\right)\dot
  S\left(\vartheta _v,X_t\right)^* }{\sigma
  \left(X_t\right)^2}\,{\rm d}t\,{\rm d}v\right|\nonumber\\
&\qquad \leq C\left|
T^{\frac{\delta }{2}} \left( \bar\vartheta_{\tau_\delta
}-\vartheta _0\right) \frac{1}{T^{\delta }}\int_{0 }^{1}\int_{0 }^{T^\delta}\frac{\dot
  S\left(\bar\vartheta_{\tau_\delta } ,X_t\right)\dot
  S\left(\vartheta _v ,X_t\right)^* }{\sigma
  \left(X_t\right)^2}\,{\rm d}t{\rm d}v\right|
T^{-\frac{1-\delta }{2}}\nonumber\\
&\qquad =O\left(T^{-\frac{1-\delta }{2}} \right)\longrightarrow 0.
\label{2-16}
\end{align}
Note that the convergence \eqref{2-15} is uniform w.r.t. $\tau \in \left[\tau
  _* ,1\right]$ because by BDG inequality: for any $\nu >0$ and any $m>0$
\begin{align*}
&\Pb_{\vartheta _0}\left(\sup_{\tau _*\leq \tau \leq 1}\frac{1}{\sqrt{T}} \left|\int_{T^\delta
    }^{\tau T}\frac{\dot S\left(\vartheta_0 ,X_t\right)
      -\dot S\left(\bar\vartheta_{\tau_\delta } ,X_t\right) }{\sigma
      \left(X_t\right)}\,{\rm d}W_t\right|>\nu  \right)\\
&\qquad \leq \frac{K_m}{\nu ^{2m}T}\int_{T^\delta
    }^{ T} \Ex_{\vartheta _0} \left|\frac{\dot S\left(\vartheta_0 ,X_t\right)
      -\dot S\left(\bar\vartheta_{\tau_\delta } ,X_t\right) }{\sigma
      \left(X_t\right)} \right|^{2m}{\rm d}t\rightarrow 0
\end{align*}
with some constant $K_m>0$.

The convergence  \eqref{2-11} follows from  the central limit theorem for the
vector-stochastic integrals (see, e.g., \cite{Kut04}, Proposition 1.21).
Further for $\tau _1<\tau _2$ we have
\begin{align*}
&\Ex_{\vartheta _0}\left|\tilde\eta _{\tau_1
  ,T}\left(\vartheta _0\right)-\tilde\eta _{\tau_2 ,T}\left(\vartheta _0\right)
\right|^4  =\frac{1}{{T^2}} \Ex_{\vartheta _0}\left|  \int_{\tau_1 T
    }^{\tau_2 T}\frac{\dot S\left(\vartheta_0 ,X_t\right)
       }{\sigma
      \left(X_t\right)}\,{\rm d}W_t \right|^4\\
&\qquad \leq \frac{\tau _2-\tau _1}{{T}} \Ex_{\vartheta _0} \int_{\tau_1 T
    }^{\tau_2 T}\left| \frac{\dot S\left(\vartheta_0 ,X_t\right)
       }{\sigma      \left(X_t\right)}\right|^4\,{\rm d}t \leq  C\left|\tau _1-\tau _2\right|^2
\end{align*} 
where the constant $C>0$ does not depend on $\vartheta _0$ and $T$.

More detailed analysis based on the same estimates shows that we have the
convergence of moments uniform on compacts 
\begin{align*}
\sup_{\vartheta _0\in {\bf A}}\Ex_{\vartheta _0}\left|\sqrt{\tau T}\II\left(\vartheta
_0\right)^{1/2} \left(\vartheta _\tau ^\star-\vartheta
_0\right)\right|^p\longrightarrow \Ex \left|\zeta \right|^p.
\end{align*}
From this convergence and the continuity of the matrix $\II\left(\vartheta
\right)$ follows the asymptotic efficiency  \eqref{2-2} of the one-step MLE.

\subsection{Two-step MLE ($\delta \in (\frac{1}{4},\frac{1}{2}]$)}

The learning time interval can be shorter. Let us take the {\it first}  estimator
$\tilde\vartheta _{\tau_\delta }$ constructed by the observations $X^{T^\delta }=\left(X_t,
,0\leq t\leq T^\delta\right) $ with $\delta
\in(\frac{1}{4},\frac{1}{2}]$. We suppose that this estimator is
consistent, asymptotically normal and the moments
converge too:
 $$
\tilde v_{\tau_\delta }=T^{\frac{\delta }{2}} \left(\tilde\vartheta
_{\tau_\delta }-\vartheta _0\right)\Longrightarrow  {\cal N} \left(0,
\MM\left(\vartheta _0\right)\right),\qquad  \sup_{\vartheta _0\in
  {\bf A}}\Ex_{\vartheta _0}\left| \tilde v_{\tau_\delta }\right|^p\leq C,
$$ 
for any $p>0$. Here $\MM\left(\vartheta _0\right) $ is some matrix and $C>0$
does not depend on $T$.  As before it can be the MLE, MDE, BE or the EMM. 

 Introduce the {\it second} preliminary estimator, which is estimator-process
\begin{equation}
\label{2-d1}
\bar\vartheta _{\tau} =\tilde\vartheta _{\tau_\delta }+\left({\tau T}\right)^{-1/2}{\II\left(
  \tilde\vartheta _{\tau_\delta }\right)^{-1}}{  }\Delta
  _{\tau T}
\left(\tilde\vartheta _{\tau_\delta },X^{\tau T}_{T^\delta
} \right),\quad \tau \in
\left[\tau _\delta ,1\right]  ,
\end{equation}
where 
\begin{align}
\Delta _{\tau T}
\left(\vartheta ,X^{\tau T}_{T^\delta
} \right)&=\frac{1}{\sqrt{{\tau T}}}\int_{T^\delta
}^{\tau T}\frac{\dot S\left(\vartheta ,X_t\right)}{\sigma
  \left(X_t\right)^2}\left[{\rm d}X_t-S\left(\vartheta ,X_t\right){\rm
    d}t\right]. 
\label{2-8a}
\end{align}
The two-step MLE-process we define as follows
\begin{equation}
\label{2-2s}
\vartheta _{\tau }^{\star\star} = \bar\vartheta _{\tau  }+\frac{{\II\left( \bar\vartheta_{\tau  }\right)^{-1}}}{\sqrt{\tau T}  }{ 
}\hat\Delta _{\tau T} \left(\tilde\vartheta _{\tau_\delta }, \bar\vartheta _{\tau  },X^{\tau T}_{T^\delta }
\right),\quad \tau _\delta \leq \tau \leq 1,
\end{equation}
where 
\begin{align*}
\hat\Delta _{\tau T}
\left(\vartheta_1 ,\vartheta_2 ,X^{\tau T}_{T^\delta } \right)&=\frac{1}{\sqrt{{\tau T}}}\int_{T^\delta
}^{\tau T}\frac{\dot S\left(\vartheta_1 ,X_t\right)}{\sigma
  \left(X_t\right)^2}\left[{\rm d}X_t-S\left(\vartheta_2 ,X_t\right){\rm
    d}t\right].
\end{align*}
Note that $\hat\Delta _{\tau T}
\left(\vartheta ,\vartheta ,X^{\tau T}_{T^\delta } \right)=\Delta _{\tau T}
\left(\vartheta ,X^{\tau T}_{T^\delta } \right) $.
\begin{theorem}
\label{T2} Suppose that the conditions of regularity hold.
 Then the two-step MLE-process  $\vartheta _{\tau }^{\star\star}, \tau _\delta \leq
\tau \leq 1$ is uniformly consistent, asymptotically normal
\begin{align*}
\sqrt{T}\left( \vartheta _{\tau }^{\star\star}-\vartheta _0 \right)\Longrightarrow
     {\cal N}\left(0,\tau ^{-1}\II\left(\vartheta _0\right)^{-1}\right),
\end{align*}
and asymptotically efficient. 
The random process
\begin{align*}
\eta _{\tau ,T}\left(\vartheta _0\right)=\tau \sqrt{T}\II\left(\vartheta
_0\right)^{-1/2}\left(\vartheta _{\tau 
}^{\star\star}-\vartheta _0  \right) ,\qquad \tau _*\leq \tau \leq 1
\end{align*}
for any $\tau _*\in \left(0,1\right)$ converges in distribution to the
$d$-dimensional standard Wiener
process $W\left(\tau \right),\tau _*\leq \tau \leq 1$. 
\end{theorem}
{\bf Proof.} The proof will be given in two steps. First we show that the
estimator-process  $\bar\vartheta _{\tau } $ is such that 
\begin{align*}
\sup_{\vartheta _0\in {\bf A}}\Ex_{\vartheta _0}\left|T^{\frac{\gamma }{2}}
\left(\bar\vartheta _{\tau }-\vartheta _0\right)\right|^p\leq C 
\end{align*}
with $\gamma \in (\frac{1}{2},1)$ and then we can use the proof of the Theorem
\ref{T1}, where the mentioned properties are already established.

 Let us take  such $\gamma >\frac{1}{2}$ that $\gamma <2\delta $.  We have
\begin{align*}
&T^{\frac{\gamma }{2}}\left(\bar\vartheta _{\tau } -\vartheta
  _0\right)=T^{\frac{\gamma }{2}}\left(\tilde\vartheta _{ \tau_\delta }
  -\vartheta _0\right)+\frac{\II\left( \tilde\vartheta
    _{\tau_\delta }\right)^{-1}}{ {\tau \sqrt{T}} }\int_{T^\delta 
  }^{\tau T}\frac{\dot S\left(\tilde\vartheta _{ \tau_\delta }
    ,X_t\right)}{\sigma \left(X_t\right)}{\rm d}W_t\;T^{\frac{\gamma-1 }{2}}\\ 
&\qquad    +T^{\frac{\gamma }{2}}\frac{\II\left( \tilde\vartheta _{\tau_\delta
    }\right)^{-1}}{ {\tau T} }\int_{T^\delta  }^{\tau
    T}\frac{\dot S\left(\tilde\vartheta _{ \tau_\delta } ,X_t\right)}{\sigma
    \left(X_t\right)^2}\left[S\left(\vartheta_0 ,X_t\right)-S(\tilde\vartheta
    _{ \tau_\delta } ,X_t)\right]{\rm d}t\\
&\quad =\hat v_{\tau_\delta }\II\left( \tilde\vartheta
    _{\tau_\delta }\right)^{-1} T^{\frac{\gamma-\delta  }{2}}\left[\II\left( \tilde\vartheta
    _{\tau_\delta }\right)-  \int_{0}^{1} \int_{T^\delta 
}^{\tau T}\frac{\dot S\left(\tilde\vartheta _{ \tau_\delta } ,X_t\right)\dot
      S\left(\vartheta_v,X_t\right)^*}{\tau T\;\sigma  
  \left(X_t\right)^2}{\rm d}t {\rm d}v \right]\\
&\qquad\quad  +O\left(T^{\frac{\gamma-1 }{2}} \right).
\end{align*}
Here  $\hat v_{\tau_\delta }= T^{\frac{\delta  }{2}}\left(\tilde\vartheta _{ \tau_\delta }
  -\vartheta _0\right)  $. We can write 
\begin{align}
&T^{\frac{\gamma -\delta }{2}}\left[\II\left( \tilde\vartheta
    _{\tau_\delta }\right)-   \int_{0}^{1}\int_{T^\delta 
}^{\tau T}\frac{\dot S\left(\tilde\vartheta _{ \tau_\delta } ,X_t\right)\dot
      S\left(\vartheta _{ v } ,X_t\right)^*}{\tau T\;\sigma 
  \left(X_t\right)^2}{\rm d}t \,{\rm d}v \right]\nonumber\\
&\qquad =\sqrt{T}\left[\II\left( \vartheta
    _0\right)-   \int_{0
}^{\tau T}\frac{\dot S\left(\vartheta _{ 0 } ,X_t\right)\dot
      S\left(\vartheta _{ 0 } ,X_t\right)^*}{\tau T\;\sigma 
  \left(X_t\right)^2}{\rm d}t  \right]T^{-\frac{1-\gamma+\delta }{2}}\nonumber\\
&\qquad\quad +\left[\II( \tilde\vartheta
    _{\tau_\delta })-\II\left( \vartheta
    _{0 }\right)\right]T^{\frac{\gamma -\delta }{2}}+\int_{0
}^{ T^\delta }\frac{\dot S\left(\vartheta _{ 0 } ,X_t\right)\dot
      S\left(\vartheta _{ 0 } ,X_t\right)^* }{\tau T\;\sigma  
  \left(X_t\right)^2}{\rm d}t\; T^{\frac{\gamma -\delta }{2}} \nonumber \\
&\qquad \quad
 +\int_{0}^{1}\int_{T^\delta 
}^{\tau T}\frac{\dot S\left(\vartheta _{ 0 } ,X_t\right)\dot
      S\left(\vartheta _{ 0 } ,X_t\right)^*-\dot S\left(\tilde\vartheta _{
        \tau_\delta  } ,X_t\right)\dot 
      S\left(\vartheta _{ v } ,X_t\right)^* }{\tau T\;\sigma  
  \left(X_t\right)^2}{\rm d}t\;{\rm d}v T^{\frac{\gamma -\delta }{2}}\nonumber\\
&\qquad =O\left( T^{-\frac{1+\gamma -\delta
    }{2}}\right)+O\left(T^{\frac{\gamma -2\delta }{2}}\right)+O\left(T^{-1+\frac{\gamma +\delta
    }{2}}  \right)+O\left( T^{\frac{\gamma -2\delta }{2}}\right)
\label{2-20}
\end{align}
Recall that the components of the vector $\dot S\left(\vartheta ,x\right)$, of
the matrix $\ddot
\SS\left(\vartheta ,x\right)$ and of the function $\sigma \left(x\right)^{-1} $ have polynomial
majorants and the invariant density has exponentially decreasing
tails. Therefore it can be shown that the moments converge too. Moreover for any $p>0$
\begin{align*}
\sup_{\vartheta _0\in {\bf A}}\Ex_{\vartheta _0}\left|T^{\gamma /2}\left(\bar\vartheta _{\tau } -\vartheta
  _0\right)\right|^p\rightarrow 0.
\end{align*}
We have the similar relations for the two-step MLE-process too. Indeed
\begin{align*}
&\sqrt{T}\left(\vartheta _{\tau }^{\star \star} -\vartheta
  _0\right)=\sqrt{T}\left(\bar\vartheta _{ \tau } -\vartheta
  _0\right)+\frac{\II\left( \bar\vartheta _{ \tau }\right)^{-1}}{ {\tau
      \sqrt{T}} }\int_{T^\delta }^{\tau T}\frac{\dot S\left(\tilde\vartheta
    _{\tau_\delta } ,X_t\right)}{\sigma \left(X_t\right)}{\rm
    d}W_t\\ 
&\qquad\quad +\frac{\II\left( \bar\vartheta
    _{\tau  }\right)^{-1}}{ {\tau \sqrt{T}} }\int_{T^\delta }^{\tau T}\frac{\dot
    S\left(\tilde\vartheta _{ \tau_\delta } ,X_t\right)}{\sigma
    \left(X_t\right)^2}\left[S\left(\vartheta_0 ,X_t\right)-S(\bar\vartheta
    _{ \tau  } ,X_t)\right]{\rm d}t\\
 &\qquad = v_{\tau 
  }^\star\II\left( \bar\vartheta _{\tau  }\right)^{-1}\left[\II\left(
    \bar\vartheta _{\tau  }\right)-\int_{0}^{1} \int_{T^\delta }^{\tau T}\frac{\dot
      S\left(\tilde\vartheta _{ \tau_\delta } ,X_t\right)\dot
      S\left(\vartheta _{ v} ,X_t\right)^*}{\tau T\;\sigma
      \left(X_t\right)^2}{\rm d}t\,{\rm d}v \right]
\\ &\qquad\quad
  +\frac{\II\left( \bar\vartheta _{ \tau }\right)^{-1}}{ {\tau
      \sqrt{T}} }\int_{T^\delta }^{\tau T}\frac{\dot S\left(\tilde\vartheta
    _{\tau_\delta } ,X_t\right)}{\sigma \left(X_t\right)}{\rm
    d}W_t,
\end{align*}
where $ v_{\tau }^\star= T^{\frac{\gamma }{2}}\left(\bar\vartheta _{ \tau }
-\vartheta _0\right) T^{\frac{1-\gamma }{2}}$. Then the corresponding
relations are
\begin{align}
&\left[\II\left( \tilde\vartheta
    _{\tau_\delta }\right)-  \int_{0}^{1} \int_{T^\delta 
}^{\tau T}\frac{\dot S\left(\tilde\vartheta _{ \tau_\delta } ,X_t\right)\dot
      S\left(\vartheta _{ v} ,X_t\right)^*}{\tau T\;\sigma 
  \left(X_t\right)^2}{\rm d}t \,{\rm d}v \right]T^{\frac{1-\gamma  }{2}}\nonumber\\
&\qquad =\sqrt{T}\left[\II\left( \vartheta
    _0\right)-   \int_{0
}^{\tau T}\frac{\dot S\left(\vartheta _{ 0 } ,X_t\right)\dot
      S\left(\vartheta _{ 0 } ,X_t\right)^*}{\tau T\;\sigma 
  \left(X_t\right)^2}{\rm d}t  \right]T^{-\frac{\gamma }{2}}\nonumber\\
&\qquad\quad +\left[\II( \bar\vartheta
    _{\tau  })-\II\left( \vartheta
    _{0 }\right)\right]T^{\frac{1-\gamma  }{2}}+\int_{0
}^{ T^\delta }\frac{\dot S\left(\vartheta _{ 0 } ,X_t\right)\dot
      S\left(\vartheta _{ 0 } ,X_t\right)^* }{\tau T\;\sigma  
  \left(X_t\right)^2}{\rm d}t\; T^{\frac{1-\gamma }{2}} \nonumber \\
&\qquad \quad +\int_{0}^{1}\int_{T^\delta 
}^{\tau T}\frac{\dot S\left(\vartheta _{ 0 } ,X_t\right)\dot
      S\left(\vartheta _{ 0 } ,X_t\right)^*-\dot S\left(\tilde\vartheta _{
        \tau_\delta  } ,X_t\right)\dot 
      S\left(\vartheta _{ v } ,X_t\right)^* }{\tau T\;\sigma  
  \left(X_t\right)^2}{\rm d}t\;{\rm d}v T^{\frac{1-\gamma  }{2}}\nonumber\\
&\qquad =O\left( T^{-\frac{\gamma 
    }{2}}\right)+O\left(T^{-\frac{2\gamma -1 }{2}}\right)+O\left(T^{-\frac{1-2\delta+\gamma 
    }{2}}  \right)+O\left( T^{-\frac{\gamma +\delta-1 }{2}}\right).
\label{2-21}
\end{align}
Therefore
\begin{align*}
\tau \sqrt{T}\left(\vartheta _{\tau }^{\star \star} -\vartheta
_0\right)&=\frac{\II\left( \vartheta _0\right)^{-1}}{ { \sqrt{T}}
}\int_{0}^{\tau T}\frac{\dot S\left(\vartheta _{0} ,X_t\right)}{\sigma
  \left(X_t\right)}{\rm d}W_t+o\left(1\right)\nonumber\\
&\Longrightarrow {\cal
  N}\left(0,\tau\II\left( \vartheta _0\right)^{-1} \right).
\end{align*}

The weak convergence $\eta _{\tau ,T}\left(\vartheta _0\right),\tau _*\leq
\tau \leq 1$ now follows from the proof of the   Theorem \ref{T1}.

\subsection{Example}

Suppose that the observed process is 
\begin{align*}
{\rm d}X_t=-\left(X_t-\vartheta \right)^3{\rm d}t+\,  {\rm d}W_t,\quad
X_0,\quad 0\leq t\leq T 
\end{align*}
where $\vartheta \in \Theta=\left(a ,b \right) $. It is easy to see
that the conditions of regularity are fulfilled and the process is ergodic
with the density of invariant distribution
\begin{align*}
f\left(\vartheta ,x\right)=\frac{8^{1/4}}{\Gamma
  \left(\frac{1}{4}\right)}\exp\left\{-\frac{\left(x-\vartheta 
  \right)^4}{2}\right\}=\varphi \left(x-\vartheta \right) .
\end{align*}
Note that the MLE of the parameter $\vartheta $ can not be written in
explicit form. Let us take $\delta =\frac{3}{4}$. We have for the empirical
mean (estimator of the method of moments) the consistency
\begin{align*}
\bar \vartheta _{T^{{3}/{4}} }=\frac{1}{T^{{{3}/{4}}}}\int_{0}^{T^{{{3}/{4}}}} X_t\;{\rm
  d}t\longrightarrow \Ex_{\vartheta _0}\xi  =\vartheta _0
\end{align*}
and asymptotic normality
\begin{align*}
T^{\frac{3}{8}}(\bar \vartheta _{T^{\frac{3}{4}}}-\vartheta _0)
=\frac{1}{T^{\frac{3}{8}} }\int_{0}^{T^{{3}/{4}} }\left( X_t-\vartheta _0\right)\;{\rm 
  d}t\Longrightarrow
 {\cal N}\left(0,D^2\right),
\end{align*}
where 
\begin{align*}
D^2=4\Ex_{\vartheta _0}\left( \int_{-\infty }^{\xi }\frac{\left(y-\vartheta
  _0\right)f\left(\vartheta _0,y \right)}{ f\left(\vartheta _0,\xi
  \right)}\;{\rm d}y\right)^2 =4\Ex_{0}\left( \int_{-\infty }^{\xi_0
}\frac{y\varphi \left(y \right)}{ \varphi \left(\xi_0 \right)}\;{\rm
  d}y\right)^2 .
\end{align*}
Here the random variable $\xi_0 $ has the density function $\varphi
\left(x\right)$. The Fisher information ${\rm I}$ does not depend on
$\vartheta $ and the one-step MLE-process is
\begin{align*}
\vartheta _\tau ^\star=\bar \vartheta _{T^{{3}/{4}} } -\frac{3}{\tau T\sqrt{{\rm I}}}\int_{T^\delta
}^{\tau T}\left( X_t-\bar \vartheta _{T^{{3}/{4}} }\right)^2\left[{\rm d}X_t+\left(
  X_t-\bar \vartheta _{T^{{3}/{4}} } \right)^3{\rm d}t\right]. 
\end{align*}
This estimator by Theorem \ref{T1} is uniformly consistent, asymptotically
normal 
\begin{align*}
\sqrt{\tau T}\left(\vartheta _\tau ^\star-\vartheta _0 \right)\Longrightarrow {\cal
N}\left(0,{\rm I}^{-1}\right)
\end{align*}
and asymptotically efficient.

If the learning interval is $\left[0,T^{3/8}\right]$, then the preliminary
estimator $\tilde\vartheta _{T^{3/8}}$ has the rate of convergence $T^{3/16}
$. We take the second estimator-process as 
\begin{align*}
\bar\vartheta _\tau =\tilde\vartheta _{T^{3/8}}+\frac{3}{{\rm I}  \tau
  T}\int_{T^{3/8}}^{\tau T}\left(X_s-\tilde\vartheta
_{T^{3/8}}\right)^2\left[{\rm d}X_s+ \left(X_s-\tilde\vartheta
_{T^{3/8}}\right)^3{\rm d}s\right] .
\end{align*}
 For this estimator the relation 
\begin{align*}
\Ex_{\vartheta _0}\left|T^{5/16}\left(\bar\vartheta _\tau -\vartheta
_0\right)\right|^p\rightarrow 0 
\end{align*} 
holds. 
Therefore by Theorem \ref{T2} the 
two-step MLE-process 
\begin{align*}
\vartheta _\tau^{\star\star} =\bar\vartheta _{\tau }+\frac{3}{{\rm I}  \tau
  T}\int_{T^{3/8}}^{\tau T}\left(X_s-\tilde\vartheta
_{T^{3/8}}\right)^2\left[{\rm d}X_s+ \left(X_s-\bar\vartheta
_{\tau }\right)^3{\rm d}s\right] 
\end{align*}
is asymptotically normal 
\begin{align*}
\sqrt{\tau T}\left(\vartheta _\tau^{\star\star} -\vartheta _0
\right)\Longrightarrow {\cal N}\left(0,{\rm I}^{-1}\right).
\end{align*}

The similar estimator-processes can be constructed and in the case of two-dimensional parameter
$\vartheta =\left(\alpha ,\beta \right)$ and the observations
\begin{align*}
{\rm d}X_t=-\beta \left(X_t-\alpha  \right)^3\,{\rm d}t+\,  {\rm d}W_t,\quad X_0,\quad 0\leq t\leq T
\end{align*}
where $\beta >0$. Indeed suppose that $\delta =\frac{3}{4}$. the preliminary
estimator $\bar\vartheta _{T^{3/4}}=\left(\bar\alpha _{T^{3/4} },\bar\beta  _{T^{3/4} } \right)$ can be
\begin{align*}
\bar\alpha _{T^{3/4}}&=\frac{1}{T^{3/4} }\int_{0}^{T^{3/4} }X_t\,{\rm
  d}t\rightarrow \alpha ,\\
\bar\beta _{T^{3/4} }&=\left(\frac{\Gamma \left(\frac{3}{4}\right)}{\Gamma
  \left(\frac{1}{4}\right)}\right)^2 \left(\frac{1}{2T^{3/4}
}\int_{0}^{T^{3/4} } \left(X_t-\bar\alpha _{T^{3/4}}\right)^2{\rm
  d}t\right)^{-2}\longrightarrow \beta .
\end{align*}
The invariant density is 
\begin{align*}
f\left(\vartheta ,x\right)=\frac{\left(8\beta \right)^{1/4}}{\Gamma
  \left(\frac{1}{4}\right)} \exp\left\{-\frac{\beta }{2}\left(x-\alpha \right)^4\right\}
\end{align*}
and the Fisher matrix $\II\left(\beta \right)$ is diagonal. Therefore the
one and two-step MLE-processes can easily by written. 

\subsection{Discussions}

Note that the process of construction of multi-step estimators can be
continued. For example, if the initial rate is $T^\delta $ with $\delta \in
(\frac{1}{8},\frac{1}{4}]$, then we can use once more one-step device to
  improve the rate of preliminary estimator up to $\gamma\in
  (\frac{1}{4},\frac{1}{2}) $, where $\gamma $ satisfies the condition
  $\gamma<2\delta $ and so on. Therefore the asymptotically efficient
  estimator-process will be three-step MLE.

We used two estimators $\tilde\vartheta _{\tau_\delta }$ and $\bar\vartheta _\tau
$ because the estimator  $\bar\vartheta _\tau
$ depends on the whole trajectory $X^{\tau T}$ and the stochastic integral
\begin{align*}
\int_{T^\delta }^{\tau T}\frac{\dot S\left(\bar\vartheta
  _\tau,X_t\right)}{\sigma \left(X_t\right)^2} \;{\rm d}X_t
\end{align*}
is not well defined. Another possibility 
to avoid this problem is to  replace the stochastic
integral by ordinary integrals and to use one estimator only as follows. 

Suppose that the functions $\dot
S\left(\vartheta ,x\right),\ddot S\left(\vartheta ,x\right)$ and $\sigma
\left(x\right)$ are continuously differentiable w.r.t. $x$ and the derivatives
belong to the class ${\cal P}$. Introduce the vector-process
\begin{align*}
\Delta^\circ _{\tau T }
\left(\vartheta,X^{\tau T}_{T^\delta} \right)&=\frac{1}{\sqrt{{\tau T}}}\int_{X_{T^\delta}
}^{X_{\tau T }}\frac{\dot S\left(\vartheta ,y\right)}{\sigma
  \left(y\right)^2}{\rm d}y -\int_{T^\delta }^{\tau T} \frac{\dot
  S'\left(\vartheta ,X_t\right) }{2\sqrt{{\tau T}}}\;{\rm d}t\\ 
& +\int_{T^\delta }^{\tau T} \frac{\dot
  S\left(\vartheta ,X_t\right)\sigma 
  \left(X_t\right)\sigma' 
  \left(X_t\right)-\dot S\left(\vartheta ,X_t\right)S\left(\vartheta
  ,X_t\right)}{\sqrt{{\tau T}}\sigma 
  \left(X_t\right)^2}\;{\rm d}t. 
\end{align*}
We have  (by the It\^o formula)
\begin{align*}
\Delta^\circ _{\tau T }
\left(\vartheta,X^{\tau T}_{T^\delta } \right)&=\Delta _{\tau T }
\left(\vartheta,X^{\tau T}_{T^\delta } \right).
\end{align*}

The two-step MLE-process in this case is 
\begin{align*}
\vartheta _{\tau }^{\circ} = \bar\vartheta _{\tau  }+\left({\tau
  T}\right)^{-1/2}{{\II\left( \bar\vartheta_{\tau  }\right)^{-1}}}{  }{  
}\Delta _{\tau T} ^\circ\left( \bar\vartheta _{\tau  },X^{\tau T}_{T^\delta }
\right),\quad \tau _\delta \leq \tau \leq 1.
\end{align*}
It can be shown that this estimator is asymptotically equivalent to $\vartheta
_{\tau }^{\star\star} $ and has the same asymptotic properties as those
described in the Theorem \ref{T2}. 

The calculation of the Fisher information matrix for some models can be a
difficult problem. In such cases we can replace the Fisher information matrix
$\II\left(\vartheta \right)$ by its empirical version
\begin{align*}
\II\left(\vartheta,t \right)=\frac{1}{t}\int_{0}^{t}\frac{\dot
  S\left(\vartheta ,X_s\right)\dot S\left(\vartheta ,X_s\right)^*}{\sigma
  \left(X_s\right)^2}\;{\rm d}s\longrightarrow \II\left(\vartheta \right) .
\end{align*}

The proposed in this work construction can be easily generalized to many other
statistical models. At particularly, it ``works'' in the case of small noise
asymptotic
\begin{align*}
{\rm d}X_t=S\left(\vartheta ,t,X_t\right)\,{\rm d}t+\varepsilon \sigma
\left(t,X_t\right)\,{\rm d}W_t,\quad X_0=x_0,\quad 0\leq t\leq T, 
\end{align*}
where $T$ is fixed and the asymptotics corresponds to $\varepsilon \rightarrow
0$. We introduce a learning time interval $\left[0,\tau _\varepsilon \right]$,
where $\tau _\varepsilon =\varepsilon ^\delta \rightarrow 0$ and for some values of $\delta
>0$ we show that the one-step MLE-process
\begin{align*}
\vartheta _{t,\varepsilon }^\star=\bar\vartheta _{\tau _\varepsilon }+\II\left(\bar\vartheta _{\tau
  _\varepsilon }\right)^{-1}\int_{\tau _\varepsilon }^{t}\frac{\dot
  S\left(\bar\vartheta _{\tau _\varepsilon },s,X_s\right)}{\sigma
  \left(s,X_s\right)^2}\left[{\rm d}X_s-S\left(\bar\vartheta _{\tau
    _\varepsilon },s,X_s\right)\right]   ,\quad t\in \left[\tau _\varepsilon
  ,T\right] 
\end{align*} 
is asymptotically efficient estimator for all $t\in \left[\tau _*,T\right]$
with any $0<\tau _*\leq T$ \cite{Kut15}. Note that for this model the
construction of the consistent preliminary estimator $\bar\vartheta _{\tau
  _\varepsilon }$  of $d$-dimensional parameter $\vartheta $ is possible if
the observed process $X_t$  is $k$ dimensional and  $k\geq d$.

This multi-step MLE-processes can be realized and in the case of estimation of
parameter $\vartheta $ by the discrete time observations
$X^n=\left(X_{t_1},\ldots,X_{t_n}\right), $ $t_j=\frac{j}{n}T$   of the diffusion 
process 
\begin{align*}
{\rm d}X_t=S\left(t,X_t\right)\,{\rm d}t+ \sigma
\left(\vartheta ,t,X_t\right)\,{\rm d}W_t,\quad X_0=x_0,\quad 0\leq t\leq T.
\end{align*}
Here we suppose that the time of observation $T$ is fixed and $n\rightarrow
\infty $. The corresponding multi-step pseudo MLE-process is  asymptotically efficient
\cite{GK15}. 

For the nonlinear autoregresive model 
\begin{align*}
X_{j+1}=S\left(\vartheta ,X_j\right)+\varepsilon _{j+1},\quad j=0,1,\ldots,n-1
\end{align*} 
the similar multi-step MLE-process provides asymptotically efficient estimator
process too \cite{KM15}.

The  construction of the  multi-step MLE-processes can be done in the case
of inhomogeneous Poisson processes, i.i.d. observations and so on.

In the work \cite{Kut15} we apply the one-step MLE  in the construction of the 
goodness-of-fit tests based on score-function-processes.

Note as well that the one-step MLE-process device allowed us to construct
asymptotically efficient estimator of the paramezters of hidden telegraph
signal \cite{KhK15}.

{\bf Aknowlledgement.} This work was done under partial financial support of
the grant of  RSF number 14-49-00079.


\begin{thebibliography}{99}
\bibitem {GK15} Gasparyan, S. and  Kutoyants, Y.A. (2015) Appproximation of
  the solution of BSDE by high frequency data. Submitted. 

 \bibitem {H67}  Holevo, A. (1967) Estimates of parameters of a diffusion
   process via stochastic approximation  
   method. {\it Proc. Comp. Center USSR Acad. Sci.},  v.12, 179-120   (rus)



\bibitem {IH81} Ibragimov I.A. and Has'minskii R.Z. (1981) {\it Statistical
Estimation - Asymptotic Theory.} {Springer-Verlag}, New York.

\bibitem{KU14} Kamatani, K. and Uchida, M. (2014) Hybrid multi-step estimators
  for stochastic differential equations based on sampled data. To appear in
  {\it Statist. Inference Stoch. Processes}


\bibitem{KS} Karatzas, I. and Shreve, S.E. (1991) {\it Brownian Motion and
  Stochastic Calculus}. 2-nd Ed., Springer, N.Y.

\bibitem{Kes97} Kessler, M. (1997) Estimating of ergodic diffusion from
  discrete observations. {\it Scand. J. Stat.}, 24, 211-229.

\bibitem{Kh12} Khasminskii, R. (2012) {\sl Stochastic Stability of
Differential Equations.}  2-nd Ed., Springer, Berlin.

\bibitem{KhK15} Khasminskii, R. and Kutoyants Yu.A. (2015) On parameter
  estimation of hidden telegraph signal. Submitted. 


\bibitem{KY} Kushner, H. and Lin, G. (2003) {\it Stochastic Approximation and
  Recursive  Algorithms and Applications.} Springer, N.Y.

\bibitem{Kut77} Kutoyants Yu.A.  (1977) { Estimation of the trend parameter of
  a diffusion process.} {\it Theory Probab. Appl}., 22, 399-405.

\bibitem {Kut04} Kutoyants, Y.A. (2004) {\it Statistical Inference for Ergodic
  Diffusion Processes.}  { Springer}, London.


\bibitem {Kut14} Kutoyants, Y.A. (2014) On approximation of the backward
  stochastic differential equation. Small noise, large samples and high
  frequency cases. {\it Proceed. Steklov Inst. Mathematics}, v. 287, 133-154.

\bibitem {Kut15} Kutoyants, Y.A. (2015) On score-functions and goodness-of-fit
  tests for stochastic processes. Submitted.

\bibitem {KM15} Kutoyants, Y.A. and Motrunich, A. (2015) On milti-step
  MLE-process for  Markov sequences.  Submitted.

\bibitem {KZ14} Kutoyants, Y.A. and Zhou, L. (2014) On approximation of the
  backward stochastic differential equation.   {\it J. Stat. Plann. Infer. }
  150, 111-123. 

\bibitem{LC56} Le Cam, L. (1956) On the asymptotic theory of estimation and
  testing hypotheses. {\it Proc. 3rd Berkeley Symposium I}, 355-368.

\bibitem{LR05} Lehmann, E.L. and Romano, J.P. (2005) {\sl Testing Statistical
Hypotheses.} (3rd ed.) Springer, N.Y.

 \bibitem {LSZ94} Levanony, D., Shwartz, A. and Zeitouni, O. (1994)
Recursive identification in continuous-time stochastic process. {\it
Stochastic Process. Appl.}, 49, 245-275.

\bibitem{LS05} Liptser, R. and Shiryaev, A.N. (2005) {\it Statistics of Random
  Processes.} v. 2, 2-nd ed. Springer, N.Y.

%\bibitem{MY} Ma, J. and Yong, J. (1999) {\it Forward-Backward Stochastic
%%  Differential Equations and their Applications.} Lecture Notes in
 % Mathematics. {Springer}, Berlin.

\bibitem{Man68}  Mandl, P.  (1968)   {\sl Analytical Treatment of One-Dimensional
 Markov Processes}. Academia, Prague; New York: Springer-Verlag.



\bibitem{NH73} Nevelson, M.B. and Hasminskii, R.Z. (1973) {\it Stochastic Approximation
and Recursive Estimation.} AMS Providence, Rhode Island.


%\bibitem {PP92} Pardoux, E. and Peng, S. (1992) Backward stochastic differential
%  equation and quasilinear parabolic differential equations. In {\it
%    Stochastic Partial Differential Equations and Their Applications, }
%  Lecture Notes in Control and Information Sciences, 176, 200-217.

\bibitem {UY12} Uchida, M. and Yoshida, N. (2012) Adaptive 
  estimation of ergodic diffusion process based on sampled data. {\it
   Stoch.  Proces. Appl.}, 122, 2885-2924. 

\bibitem {UY14} Uchida, M. and Yoshida, N. (2014) Adaptive Bayes type
  estimators of ergodic diffusion processes from discrete observations. {\it
   Statist. Inference  Stoch. Processes.}, 17, 3, 181-219. 

\bibitem {Y92} Yoshida, N. (1992) Estimation for diffusion processes from
  discrete observations. {\it J. Multivariate Anal.}, 41, 220-242.

\end{thebibliography}
\end{document}